\UseRawInputEncoding %arxiv
\documentclass[11pt, a4paper]{article}
\usepackage[T1]{fontenc}
\usepackage[utf8]{inputenc}
\usepackage{cite}
\usepackage{amsfonts}
\usepackage{amssymb}
\usepackage{enumerate}
\usepackage{mathrsfs}
\usepackage{longtable}
\usepackage{amsmath}
\usepackage{amsthm}
\usepackage{booktabs}
\usepackage{epsf, epsfig, amsfonts, amsgen, indentfirst}
\usepackage{amsmath, amstext, amsbsy, amsopn, amsthm, bbding, wasysym}
\usepackage{multicol, mathdots}
\usepackage{tikz}
\usepackage{array}
\usepackage{setspace}% Controls spacing in selected parts

\allowdisplaybreaks
\makeatletter
\renewcommand{\@makefntext}[1]{\noindent #1}
\makeatother

\theoremstyle{plain}
\newtheorem{theorem}{Theorem}[section]

\newtheorem{lemma}{Lemma}[section]

\newtheorem{remark}{Remark}[section]

\newtheorem{conjecture}{Conjecture}[section]

\newenvironment{proofthm}[1]{\begin{proof}[\textbf{#1:}]}{\end{proof}}

\allowdisplaybreaks

\title
{\LARGE \textbf{Upper bounds for the number of spanning forests of a graph
\thanks{ supported by NSFC (No. 12261071) and NSF of Qinghai Province (No. 2025-ZJ-902T).} 
}}
\author{Tingzeng Wu$^{a,b}$\thanks{{Corresponding author.\newline
\emph{E-mail address}: mathtzwu@163.com, sxlu863@163.com, a3566293588@163.com
}}, Shuxin Lu$^a$, Xiangshuai Dong$^c$ \\
{\small $^{a}$ School of Mathematics and Statistics, Qinghai Minzu University, }\\
{\small  Xining, Qinghai 810007, P.R.~China} \\
{\small $^{b}$ Qinghai Institute of Applied Mathematics,   Xining, Qinghai, 810007, P.R.~China}\\
{\small $^{c}$ School of Mathematical Sciences, Xiamen University,  }\\
{\small  Fujian, Xiamen 361005,  P.R.~China} \\
}
\date{}
\begin{document}
\maketitle

\noindent {\bf Abstract: }
Let $G$ be a simple graph on $n$ vertices, and let $F(G)$ denote the number of spanning forests of $G$. Bencs and Csikv\'ari  proved that every $r$-regular graph $G$ with $r\ge 2$ satisfies $F(G)\le r^n$. Motivated by this result, they conjectured that if $G$ has minimum degree at least $3$,
then $F(G)\leq \prod_{i=1}^n d_i$, where $(d_1, d_2, \dots, d_n)$ is the degree sequence of $G$.
In this paper, we prove that graphs with minimum degree at least $n/3$ $(n\geq7)$
 satisfy $F(G)\leq \prod_{i=1}^n d_i$.
We also construct a counterexample to the conjecture.
We therefore further study the enumeration of spanning forests in general graphs, proving an upper bound in terms of the degree sequence
 and characterizing the extremal graphs attaining this bound.
In addition, using the matrix-forest theorem, we obtain an upper bound for the number of spanning forests of a connected graph in terms of its algebraic connectivity,
 and characterize the extremal graphs attaining this bound.
\\
\noindent {\bf Keywords: } Spanning forest; Pseudo-forest; Graph polynomial; Laplacian spectrum\\
\noindent {\bf 2020 Mathematics Subject Classification: } 05C07; 05C30; 05C35; 05C50
%05C07 Vertex degrees
%05C30 Enumeration in graph theory
%05C35 Extremal problems in graph theory
%05C50 Graphs and linear algebra (matrices, eigenvalues, etc.)
\section{Introduction}
Throughout this paper, we consider only finite, simple, undirected graphs. Let $G = (V(G), \\ \ E(G))$ be a graph with vertex set $V(G)=\{ v_1, v_2, \dots, v_n \}$ and edge set $E(G)=\{ e_1, e_2, \dots, e_m \}$. The numbers of vertices and edges are denoted by $|V(G)| =n$ and $|E(G)| =m$, respectively. If the number of edges of $G$ is $m=0$, then $G$ is called an edgeless graph. For a vertex $v_i$, let $d_i=d_{G}(v_i)$ be its degree, that is, the number of edges incident with it. For $v \in V(G)$, let $N_G(v)$ denote the set of all vertices in $G$ adjacent to $v$. Arranging all vertex degrees in nonincreasing order gives the degree sequence of $G$, $(d_1, d_2, \dots, d_n)$, where $d_1\ge d_2\ge\cdots\ge d_n$. If $d_1 = d_2 = \cdots = d_n = r$, then $G$ is called an $r$-regular graph.
If $r= n-1$, then $G$ is called the complete graph of order $n$, denoted by $K_{n}$. Let $\delta=d_n$ denote the minimum degree of $G$. If $V(G)$ can be partitioned into two nonempty disjoint sets $V_1(G)$ and $V_2(G)$, and the edges of the graph are precisely all possible edges joining vertices in these two sets, then $G$ is called a complete bipartite graph. When $|V_1(G)|=s$ and $|V_2(G)|=t$, it is denoted by $K_{s,t}$. The graph $K_{1,t}$ $(t\geq 1)$ is a star, whose vertex of degree $t$ is its center. When $t=1$, either endpoint can serve as the center. The double star $D_{p,q}$ is obtained from two stars $K_{1,p}$ $(p\geq 1)$ and $K_{1,q}$ $(q\geq 1)$ by joining their centers with an edge. Let $G(t_1,t_2,t_3)$ denote the connected graph with $n=t_1+t_2+t_3+2$ vertices and $m=t_1+t_2+2t_3+1$ edges, obtained from the stars $K_{1,t_1+t_3}$ and $K_{1,t_2+t_3}$ by joining their centers and pairwise identifying $t_3$ leaves of $K_{1,t_1+t_3}$ with $t_3$ leaves of $K_{1,t_2+t_3}$, where $t_1,t_2\geq 0$ and $t_3\geq 1$.

For any subset $S\subseteq V(G)$, retain all vertices in $S$ and all edges of $G$ both of whose endpoints belong to $S$. The resulting subgraph is called the subgraph induced by $S$ and is denoted by $G[S]$. When $S=\varnothing$, we adopt the convention that $G[\varnothing]$ is the graph with $0$ vertices. Deleting from $G$ the vertices in $S$ and all edges incident with these vertices yields the graph denoted by $G-S$; that is, $G-S=G[V(G)\setminus S]$. In particular, when $S=\{v_i\}$, we abbreviate $G-\{v_i\}$ as $G-v_i$. For any notation not defined here, we refer the reader to \cite{Diestel2017}.

A spanning forest of $G$ is a subgraph $(V(G), E')$ with $E'\subseteq E(G)$ that is acyclic; that is, it is obtained by retaining all vertices of $G$ and choosing an edge set $E'\subseteq E(G)$ such that $(V(G), E')$ is acyclic. The set of all spanning forests of $G$ and its cardinality are denoted, respectively, by
$$
\mathcal{F}(G)=
\{(V(G), E'): E'\subseteq E(G), \ (V(G), E')\text{ is acyclic}\},
\qquad F(G)=|\mathcal{F}(G)|.
$$

Determining the number of spanning forests of a graph is a classical problem in combinatorial enumeration; see
\cite{Aldous2007, Annan1994, Bencs2022, Caracciolo2004, Grimmett2004, Liu1981, Luo2026, Teranishi2005, Thomassen2010, Wu2026, Xu2026}.
Jaeger et al. \cite{Jaeger1990} proved that computing $F(G)$ is \#P-hard.
In $2021$, Borb\'enyi et al. established the following upper bound for the number of spanning forests.

\begin{theorem}([\citen{Borbenyi2021}])
Let $G$ be a graph on $n$ vertices. Then
$$
F(G)\leq \prod_{i=1}^n (d_i+1).
$$
\end{theorem}
\noindent For regular graphs, Bencs and Csikv\'ari studied the number of spanning forests in $2023$
and obtained the following upper bound.
\begin{theorem}([\citen{Bencs2023}])\label{the:1}
Let $G$ be a graph on $n$ vertices that is $r$-regular $(r \geq 2)$. Then
$$F(G) \le r^{n}. $$
\end{theorem}
Based on the two theorems above,
 they proposed the following conjecture.
\begin{conjecture}\label{con:1}
Let $G$ be a graph on $n$ vertices. If $\delta \geq 3$,
then
$$
 F(G)\leq \prod_{i=1}^n d_i.
$$
\end{conjecture}
We investigate this conjecture and identify a class of graphs for which it holds.
\begin{theorem}\label{cor:third}
Let $G$ be a graph on $n\ge7$ vertices. If $\delta\ge n/3$, then
$F(G)<\prod_{i=1}^n d_i$.
\end{theorem}
 We also find a class of graphs that do not satisfy Conjecture \ref{con:1}.
Let $K_{y_1, y_{2}}(y_{2}> y_{1}\geq3)$ be a complete bipartite graph whose vertex set is partitioned into
 $V_1(G)=\{v_i: d_i=y_{2},1\leq i\leq y_{1}+y_{2}\}$ and $V_2(G)=\{v_j: d_j=y_{1},1\leq j\leq y_{1}+y_{2}\}$.
 We now construct a class of spanning forests of $K_{y_{1}, y_{2}}$.
For each vertex in $V_2(G)$, independently choose either to remain isolated or to be joined to one vertex in $V_1(G)$. Each such vertex has $y_1+1$ choices, and the resulting spanning subgraph is necessarily a spanning forest.
Hence
 $
F(K_{y_{1}, y_{2}})\ge(y_{1}+1)^{y_{2}},
$
and thus
$$
\frac{F(K_{y_{1}, y_{2}})}{\prod_{i=1}^{n}d_i}
\ge\frac{(y_{1}+1)^{y_{2}}}{y_{2}^{y_{1}} y_{1}^{y_{2}}}
=\frac1{y_{2}^{y_{1}}}\left(\frac{y_{1}+1}{y_{1}}\right)^{y_{2}}.
$$
Let $y_{1}$ be a fixed constant with $y_{1}\geq3$. Since, as $y_{2}\to\infty$, we have $\ln y_{2}/y_{2}\to0$, it follows that, as $y_{2}\to\infty$,
$$
\frac{1}{y_{2}^{y_{1}}}\left(\frac{y_{1}+1}{y_{1}}\right)^{y_{2}}
= \mathrm{e}^{\,y_{2}\ln\left(1+\frac{1}{y_{1}}\right)-y_{1}\ln y_{2}}
\longrightarrow +\infty.
$$
Therefore, for every fixed constant $y_{1}$, all sufficiently large integers $y_{2}$ satisfy
$$
F(K_{y_{1},y_{2}})>\prod_{i=1}^{n}d_i.
$$
As a concrete example, consider $K_{3,14}$. Applying Stark's exact
enumeration formula for spanning forests of complete bipartite graphs
\cite{Stark2013}, we obtain
$$
F(K_{3,14})
=13\,186\,891\,776
>13\,124\,466\,936
=14^3\,3^{14}.
$$

In view of these results, determining a sharp upper bound for the number of spanning forests of a general graph in terms of its degree sequence is an interesting problem.
We investigate this problem and obtain the following theorem.
\begin{theorem}\label{the:2}
Let $G$ be a graph on $n\ge1$ vertices.
For $2\le i\le n$, define
$k_i=(d_i-i+1)_+$ and $h_i=\sum_{j=i+1}^{i+k_i}\frac1{d_j+1}$,
 where $(x)_+ =\max\{x, 0\}$,
and $k_i=0$ gives $h_i=0$. Then
\begin{equation}\label{eq:1.3}
F(G)\le \prod_{i=2}^n(d_i+1-2h_i),
\end{equation}
with equality if and only if $G$ is isomorphic to an edgeless graph, a star $K_{1, t}(t\geq1)$, a double star $D_{p, q}$, a graph $G(t_1,t_2,t_3)$, or the disjoint union of one of these graphs with isolated vertices.
\end{theorem}
Chebotarev and Shamis studied the number of spanning forests of a graph using its Laplacian spectrum
 and obtained the following upper bound.
\begin{theorem}([\citen{Chebotarev1997}])\label{the:10}
Let $G$ be a graph on $n$ vertices,
 and let $L(G)$ be its Laplacian matrix,
 with eigenvalues $\lambda_1\geq\lambda_2\geq\cdots\geq\lambda_n=0$.
 Let $I_n$ denote the identity matrix of order $n$. Then
$$F(G)\leq \det(I_n+L(G))=\prod_{i=1}^{n}(1+\lambda_i),$$
 with equality if and only if $G$ is an edgeless graph.
\end{theorem}
Using Theorem \ref{the:10}, we further characterize the relationship between the number of spanning forests of a connected graph and its algebraic connectivity, obtaining the following theorem.
\begin{theorem}\label{cor:1.1}
Let $G$ be a connected graph on $n\ge2$ vertices with algebraic connectivity $\mu=\lambda_{n-1}>0$
 and $\tau(G)$ spanning trees. Then
\begin{equation}\label{eq:1.6}
F(G)\le\tau(G)((1+\mu^{-1})^n-\mu^{-n}),
\end{equation}
with equality if and only if $G\cong K_2$ or $G\cong K_3$.
\end{theorem}
The remainder of this paper is organized as follows. In Section $2$, we prove Theorem \ref{cor:1.1}.
In Section $3$, we prove Theorem \ref{cor:third}.
In Section $4$, we prove Theorem \ref{the:2}.
Finally, we further discuss the comparison between the upper bounds in Theorem \ref{the:2} and Theorem \ref{cor:1.1} for connected graphs.
\section{Proof of Theorem \ref{cor:1.1}}
To prove Theorem \ref{cor:1.1}, we first introduce a lemma.

\begin{lemma}([\citen{Kelmans1974}])\label{lem:4.1}
Let $G$ be a graph on $n$ vertices, and let $L(G)$ be its Laplacian matrix.
For $0\le k\le n-1$, write
$$
\mathcal{F}_k(G)=\{\varphi\in\mathcal{F}(G): |E(\varphi)|=k\}, \qquad f_k(G)=|\mathcal{F}_k(G)|.
$$
Let
\begin{equation}\label{eq:4.1}
\det(I_n+tL(G))=\sum_{k=0}^{n-1}c_k(G)t^k,
\end{equation}
then
\begin{equation*}
c_k(G)=\sum_{\varphi\in\mathcal{F}_k(G)}\prod_{T\in\mathcal{T}(\varphi)}|V(T)|,
\end{equation*}
where $\mathcal{T}(\varphi)$ denotes the set of all tree components of the forest $\varphi$.
 In particular, $c_0(G)=1$.
\end{lemma}
We now prove Theorem \ref{cor:1.1}.
\begin{proofthm}{Proof of Theorem \ref{cor:1.1}}
For a graph on $n$ vertices, denoted by $G$,
 let $\varphi\in\mathcal{F}_k(G)$ be a spanning forest of $G$ with $k$ edges, where $0\leq k\leq n-1$.
 Thus, $|\mathcal{T}(\varphi)|=n-k$.
Let the orders of its tree components be $s_1,\ldots,s_{n-k}$,
and let $s_j'=s_j-1\geq0$.
Since the number of edges in each tree component is $s_j-1$, we have $\sum_{j=1}^{n-k}s_j'=k$.
 Therefore,
\begin{equation}\label{eq:4.5}
 \prod_{T\in\mathcal{T}(\varphi)}|V(T)|
 =\prod_{j=1}^{n-k}(1+s_j')
 \geq1+\sum_{j=1}^{n-k}s_j'=k+1.
\end{equation}
If at least two of the $s_j'$ are positive, then the corresponding quadratic term in the expansion of the product is positive, and the inequality in \eqref{eq:4.5} is strict.
Conversely, if at most one of the $s_j'$ is positive, then all terms of degree at least $2$ are $0$.
Thus, equality holds in \eqref{eq:4.5} if and only if $\varphi$ has at most one tree component containing an edge.
By Lemma \ref{lem:4.1}, for every $0\le k\le n-1$, we have
\begin{equation}\label{eq:4.6}
c_k(G)=\sum_{\varphi\in\mathcal{F}_k(G)}\prod_{T\in\mathcal{T}(\varphi)}|V(T)|\ge(k+1)f_k(G).
\end{equation}

The Laplacian matrix is real symmetric and positive semidefinite, and $L(G)\mathbf{1}=0$,
 where $\mathbf 1$ denotes the $n\times1$ vector whose entries are all $1$.
 Therefore,
$$
\det(I_n+tL(G))=\prod_{i=1}^n(1+t\lambda_i).
$$
Combining \eqref{eq:4.1} and \eqref{eq:4.6} and integrating the preceding identity term by term yields
\begin{equation}\label{eq:1.5}
F(G)=\sum_{k=0}^{n-1}f_k(G)\le\sum_{k=0}^{n-1}\frac{c_k(G)}{k+1}
=\int_0^1\sum_{k=0}^{n-1}c_k(G)t^k\,\mathrm{d}t=
\int_0^1\prod_{i=1}^n(1+t\lambda_i)\,\mathrm{d}t,
\end{equation}
Equality holds in \eqref{eq:1.5} if and only if every spanning forest of $G$, denoted by $\varphi$, attains equality in \eqref{eq:4.5},
 that is, $\varphi$ has at most one component containing an edge,
 which is equivalent to the matching number of $G$ being at most one.

We now consider a connected graph on $n\geq2$ vertices, denoted by $G$. By Lemma \ref{lem:4.1}, we have
 $
\prod_{i=1}^{n-1}\lambda_i=c_{n-1}(G)=n\tau(G).
$ 
For $1\leq i\leq n-1$, we have $\lambda_i\ge\mu>0$. Combining this with \eqref{eq:1.5} gives
\begin{align}
F(G)
&\le n\tau(G)\int_0^1\prod_{i=1}^{n-1}(t+\lambda_i^{-1})\,\mathrm{d}t \label{eq:666}\\
&\le n\tau(G)\int_0^1(t+\mu^{-1})^{n-1}\,\mathrm{d}t \label{eq:667}\\
&= \tau(G)\bigl((1+\mu^{-1})^n-\mu^{-n}\bigr), \nonumber
\end{align}
This proves \eqref{eq:1.6}.

We now discuss the equality case in \eqref{eq:1.6},
 which is equivalent to equality holding simultaneously in \eqref{eq:666} and \eqref{eq:667}.
First, consider equality in \eqref{eq:667}.
Then $\lambda_1=\cdots=\lambda_{n-1}=\mu$,
 so $L(G)=\mu(I_n-J_n/n)$, where $J_n$ is the all-ones matrix.
 Since the off-diagonal entries of the Laplacian matrix $L(G)$ can only be $-1$ or $0$,
 we obtain $\mu=n$, that is, $G\cong K_n$.
 Equality in \eqref{eq:666} is equivalent to equality in \eqref{eq:1.5},
 that is, the maximum matching of $G$ has at most $1$ edge.
 Thus, equality holds in \eqref{eq:1.6} if and only if $G\cong K_2$ or $G\cong K_3$. This completes the proof.
\end{proofthm}

\section{Proof of Theorem \ref{cor:third}}
To prove Theorem \ref{cor:third}, we first introduce some lemmas.
\begin{lemma}([\citen{Horn2013}])\label{lem:fenkuaijuzhen}
Let $M$ be a real symmetric positive definite matrix of order $n\geq2$, 
with the following block form:
$$
M=
\begin{pmatrix}
a&\mathbf z^T\\
\mathbf z&H
\end{pmatrix},
$$
where $a$ is a real number, $\mathbf z$ is an $(n-1)\times1$ column vector,
and $H$ is a matrix of order $n-1$. Then
\begin{equation}\label{eq:fischer}
 \det M\le a\det H
 \le a\left(\frac{\operatorname{tr} H}{n-1}\right)^{n-1},
\end{equation}
Equality in the first inequality holds if and only if
 $\mathbf z=\mathbf 0$; equality in the second inequality holds if and only if
 $H$ is a positive real multiple of the identity matrix.
 \end{lemma}
\begin{lemma}\label{thm:main}
Let $G$ be a graph on $n\ge2$ vertices with no isolated vertices. Define
$$
W_n=\left(\frac{n}{n-1}\right)^{n-1},\qquad \beta(G)=\frac1n\sum_{i=1}^{n}\frac1{d_i}-\frac1{2m}.
$$
 Then
\begin{equation}\label{eq:main}
 \frac{F(G)}{\prod_{i=1}^n d_i}
 \le \frac{W_n}{2m}\left((1+\beta(G))^n-\beta(G)^n\right),
\end{equation}
Equality in~\eqref{eq:main} holds if and only if $G\cong K_2$ or $G\cong K_3$.

\end{lemma}

\begin{proof}
Let $D(G)=\operatorname{diag}(d_1, \ldots, d_n)$ denote the degree matrix of $G$, abbreviated as $D$.
 Define
$$\mathcal L=D^{-1/2}L(G)D^{-1/2}, \qquad  M_t=D^{-1}+t\mathcal L.$$
Since $\mathcal L$ is positive semidefinite and $D^{-1}$ is positive definite,
 for every real number $t\ge0$, the matrix $M_t$ is positive definite. 
Since $
D^{1/2}\mathcal LD^{1/2}=L(G)
$, we have
\begin{equation}\label{eq:juzhendengshi}
D^{1/2}M_tD^{1/2}=I_n+tL(G).
\end{equation}
Since $
\det(D)=\prod_{i=1}^n d_i,
$
taking determinants on both sides of \eqref{eq:juzhendengshi} gives
\begin{equation}\label{eq:11}
\det(M_t)=\frac{\det(I_n+tL(G))}{\prod_{i=1}^n d_i}.
\end{equation}
Let $\mathbf u=\frac{D^{1/2}\mathbf 1}{\sqrt{2m}}$ be an $n\times1$ column vector. 
Since
$$
{\mathbf u}^T \mathbf u=\frac{\mathbf 1^TD\mathbf 1}{2m}=1,
$$
$\mathbf u$ is a unit vector. Since
$$
\mathcal L \mathbf u
=\frac{D^{-1/2}L(G)D^{-1/2}
D^{1/2}\mathbf 1}{\sqrt{2m}}
=\mathbf 0,
$$
we have $M_t \mathbf u
=(D^{-1}+t\mathcal L)\mathbf u
=D^{-1}\mathbf u$.
Let $q=\mathbf u^TD^{-1} \mathbf u$.
Then$$
q
=\frac{1}{2m}
\mathbf 1^TD^{1/2}D^{-1}D^{1/2}\mathbf 1
=\frac{1}{2m}\mathbf 1^T\mathbf 1
=\frac{n}{2m}.
$$
Let $\mathbf u^\perp=\{\mathbf x\in\mathbb R^n: \mathbf u^T\mathbf x=0\}$.
Choose an $n\times(n-1)$ real matrix $Q$ whose columns form an orthonormal basis of $\mathbf u^\perp$.
Then $U=(\mathbf u, Q)$ is an orthogonal matrix. Hence
$$
 U^{\mathsf T} M_tU=
 \begin{pmatrix}
  q&\mathbf b^{\mathsf T}\\ \mathbf b&H_t
 \end{pmatrix},
 \qquad
 q=\mathbf u^{\mathsf T} D^{-1}\mathbf u=\frac{n}{2m},
 \quad \mathbf b=Q^{\mathsf T}M_t \mathbf u=Q^{\mathsf T}D^{-1}\mathbf u,
  \quad H_t=Q^{\mathsf T} M_tQ.
$$
Since $M_t$ is positive definite and the columns of $Q$ are linearly independent, $H_t$ is also positive definite.
Write $s'=\sum_i^{n} d_i^{-1}$.
Then
$$
\operatorname{tr}(M_t)
=\operatorname{tr}(D^{-1})
+t\operatorname{tr}(\mathcal L)
=s'+nt.
$$
Orthogonal similarity transformations preserve the trace of a matrix.
Combining this with $\operatorname{tr}(U^TM_tU)=q+\operatorname{tr}(H_t)$,
we obtain
\begin{equation*}
 \operatorname{tr} H_t=\operatorname{tr} M_t-q=s'+nt-\frac{n}{2m}
 =n\bigl(t+\beta(G)\bigr).
\end{equation*}
Since $H_0$ is also positive definite, $n\beta(G)=\operatorname{tr} H_0>0$.
By Lemma \ref{lem:fenkuaijuzhen}, we have
$$
 \det (U^{\mathsf T} M_tU)=
 \det M_t
 \le\frac{n}{2m}
 \left(\frac{n(t+\beta(G))}{n-1}\right)^{n-1}.
$$
Combining this with \eqref{eq:11}, we obtain, for every real number $t\ge0$,
\begin{equation}\label{eq:determinant}
 \frac{\det(I_n+tL(G))}{\prod_{i=1}^n d_i}
 \le \frac{n}{2m}
 W_n
 \bigl(t+\beta(G)\bigr)^{n-1}.
\end{equation}

From \eqref{eq:1.5} and~\eqref{eq:determinant}, we obtain
\begin{align}
\frac{F(G)}{\prod_{i=1}^n d_i}
&\le \int_0^1 \frac{\det(I_n+tL(G))}{\prod_{i=1}^n d_i}\,\mathrm{d}t \label{eq:2.9} \\
&\le \frac{n}{2m}W_n\int_0^1 (t+\beta(G))^{n-1}\,\mathrm{d}t \label{eq:2.90} \\
&= \frac{W_n}{2m}\left((1+\beta(G))^n-\beta(G)^n\right), \nonumber
\end{align}
which proves~\eqref{eq:main}.

Equality in~\eqref{eq:main} holds if and only if equality holds in both \eqref{eq:2.9} and \eqref{eq:2.90}. 
We first consider the case of equality in \eqref{eq:2.90}. 
In this case,
$$
\int_0^1 \left[ \frac{n}{2m}W_n(t+\beta(G))^{n-1}-
\frac{\det(I_n+tL(G))}{\prod_{i=1}^n d_i} \right] \,\mathrm{d}t=0.
$$
By~\eqref{eq:determinant}, the function in square brackets above
is continuous and nonnegative on $[0, 1]$. Hence, for every real number $t'\in[0, 1]$, we have
\begin{equation}\label{eq:main1}
 \frac{nW_n}{2m}(t'+\beta(G))^{n-1}=
\frac{\det(I_n+t'L(G))}{\prod_{i=1}^n d_i}.
\end{equation}
For \eqref{eq:main1} to hold, both equality conditions in \eqref{eq:fischer} must be satisfied, namely,
 $\mathbf b=Q^{\mathsf T}D^{-1}\mathbf u=\mathbf0$ and $H_t$ is a positive real multiple of the identity matrix.
 From $\mathbf b=\mathbf0$, it follows that $D^{-1}\mathbf u$ is parallel to $\mathbf u$.
 Since $G$ has no isolated vertices and every coordinate of $\mathbf u$ is nonzero,
 all $d_i$ are equal. Let $d_i=r$; that is, $G$ is $r$-regular.
 Taking $t'=1$, since $H_1$ is a positive real multiple of the identity matrix and
 $M_1=r^{-1}I_n+r^{-1}L(G)$, 
we obtain
$$
L(G)=\frac{nr}{n-1}
\left(I_n-\frac1nJ_n\right).
$$
Since every off-diagonal entry of the Laplacian matrix $L(G)$ is either $-1$ or $0$, we have $r=n-1$, 
that is, $G\cong K_{n}$.
Equality in \eqref{eq:2.9} holds if and only if equality holds in \eqref{eq:1.5}.
 Therefore, equality in \eqref{eq:main} holds if and only if $G\cong K_2$ or $G\cong K_3$.
  This completes the proof.
\end{proof}
We now prove Theorem \ref{cor:third}.
\begin{proofthm}{Proof of Theorem \ref{cor:third}}
Let $G$ be a graph on $n\ge7$ vertices with $\delta\ge n/3$. 
For a real number $x>0$, define
\begin{equation}\label{eq:B}
 B_n(x)=\frac{W_n}{nx}
 \left[
  \left(1+\frac{n-1}{nx}\right)^n
  -\left(\frac{n-1}{nx}\right)^n
 \right].
\end{equation}
For every $i(1\leq i \leq n)$, we have $0<d_i<2m$.
 Combining $2m=\sum_{i=1}^{n}d_i\geq n\delta$ with the preceding relation,
  we obtain
\begin{align}
 n\beta(G)
 =\sum_{i=1}^{n}\frac{1-d_i/(2m)}{d_i}
 \le\frac1\delta\sum_{i=1}^{n}\left(1-\frac{d_i}{2m}\right)
 =\frac{n-1}{\delta}.\label{eq:beta-delta}
\end{align}
The function $f(x)=(1+x)^n-x^n$ is strictly increasing for $x>0$, 
whereas $\frac{n-1}{nx}>0$ and is strictly decreasing for $x>0$.
Therefore,
$$
 \left(1+\frac{n-1}{nx}\right)^n
  -\left(\frac{n-1}{nx}\right)^n
$$
is strictly decreasing for $x>0$.
 Since $\frac{W_n}{nx}$ is also strictly decreasing for $x>0$, 
$B_n(x)$ is strictly decreasing for $x>0$.
 Substituting~\eqref{eq:beta-delta}
into~\eqref{eq:main} gives
\begin{equation*}
 \frac{F(G)}{\prod_{i=1}^n d_i}
 \le \frac{W_n}{2m}\left((1+\beta(G))^n-\beta(G)^n\right)
 \le B_n(\delta).
 \end{equation*}
Since $W_n<e$ and $1+x<e^x$ $(x>0)$, we have 
$$
 B_n(\delta)
 <\frac{e}{n\delta}
 \left(1+\frac{n-1}{n\delta}\right)^n
 <\frac{e^{1+(n-1)/\delta}}{n\delta}.
$$
By the above discussion, either $B_n(\delta)<1$ or $\ln(n\delta)\ge1+\frac{n-1}{\delta}$ 
 guarantees that Conjecture~\ref{con:1} holds.

We first consider the case $n\ge12$. Since $\delta\ge n/3$, we have
$$
 \ln(n\delta)-1-\frac{n-1}{\delta}
 \ge\ln(n^2/3)-4+\frac3n.
$$
Let 
$
 h(x)=\ln(x^2/3)-4+\frac3x.
$ 
 For $x\ge12$, we have $h'(x)=(2x-3)/x^2>0$, and
$h(12)=\ln48-15/4>0$, 
so, for $n\ge12$, Conjecture~\ref{con:1} holds.

When $7\le n\le11$, the integrality of the degrees
gives $\delta\ge\lceil n/3\rceil$.
Direct substitution into~\eqref{eq:B} gives (with the results rounded to two decimal places)
 $B_7(3)=0.70$, $B_8(3)=0.82$, $B_9(3)=0.98$, $B_{10}(4)=0.49$, $B_{11}(4)=0.56$.
Thus, for all $7\le n\le11$, we have $B_n(\delta)<1$, completing the proof of Theorem \ref{cor:third}.
\end{proofthm}

\section{Proof of Theorem \ref{the:2}}
To prove Theorem \ref{the:2}, we first introduce some definitions, notation, and lemmas.

A graph is called a pseudo-forest if each of its connected components is a tree or a unicyclic graph. For a graph $G$, let $\mathcal{P}(G)$ denote the set of its spanning pseudo-forests,
 let $z(\psi)$ denote the number of unicyclic components of a pseudo-forest $\psi$, and let $\mathcal{T}(\psi)$ denote the set of its tree components.
Let $M$ be a subset of $E(G)$.
If no two edges in $M$ are adjacent in $G$, then $M$ is called a matching of $G$.
Let $|M|$ denote the number of edges in $M$, and let $V(M)$ denote the set of all vertices covered by the matching $M$.
Let $\mathcal M(G)$ denote the set of all matchings of $G$ (including the empty matching).

\begin{lemma}\label{lem:2.1}
Let $G$ be a graph on $n$ vertices. For each vertex $v_i\in V(G)$,
 let $d_{i}$ denote the degree of $v_i$. Let $a_{i}\ge0$ be the weight of $v_i$, 
 and let $q_{i}=d_{i}+1+a_{i}$.
 Define
\begin{equation}\label{eq:2.1}
P_G(\mathbf{q})=\sum_{M\in\mathcal{M}(G)}(-2)^{|M|}\prod_{v_i\in V(G)\setminus V(M)}q_{i},
\end{equation}
where $\mathbf{q}$ is the parameter vector, that is, $\mathbf{q}=(q_{i})_{v_i\in V(G)}$. Then
\begin{equation}\label{eq:2.2}
P_G(\mathbf{q})=\sum_{\psi\in\mathcal{P}(G)}2^{z(\psi)}
\prod_{T\in\mathcal{T}(\psi)}\left(1+\sum_{v_i\in V(T)}a_{i}\right)>0,
\end{equation}
When $V(G)=\varnothing$, both sides are taken to be $1$.
\end{lemma}
\begin{proof}
Let $t_{i}'=1+a_{i}$. Then $q_{i}=t_{i}'+d_{i}$.
 In \eqref{eq:2.1}, expand the term corresponding to a matching $M'$ as follows: each matching edge $v_{i}v_{j}$ produces arrows $v_{i}\to v_{j}$ and $v_{j}\to v_{i}$; mark this pair of arrows and assign it weight $-2$.
 In the graph $G$, for the matching $M'$, each uncovered vertex $v_i$ has two choices:
 either choose a neighbor $v_{j}$ and draw an arrow $v_i \to v_j$, with weight $1$,
  or choose no neighbor and draw no arrow, with weight $t_{i}'$.
 The sum of the weights of all marked configurations obtained in this way is $P_G(\mathbf{q})$.
 Denote the new graph obtained from $G$ by the above procedure by $G'$.
  A vertex of $G'$ with no outgoing arrow is called a root of $G'$.
    In the graph $G'$, let the total number of pairs of opposite arrows be $b$. 
    
   Fix a graph $G'$. Each pair of opposite arrows can arise in two ways:
   either it is produced by a matching edge $v_{i}v_{j}$; 
   or, for the matching $M'$, two uncovered vertices $v_i$ and $v_{j}$ choose each other as neighbors and draw arrows $v_i \to v_j$ and $v_{j}\to v_{i}$, respectively.
 Since each vertex has at most one outgoing arrow, these pairs of opposite arrows are pairwise vertex-disjoint.
 Hence each pair can be chosen independently to arise from the matching or not.
 Let $J\subseteq\{1, \ldots, b\}$ be the set of indices of the arrow pairs marked as belonging to the matching.
  Then the weight of the part marked as belonging to the matching is
 $
(-2)^{|J|}.
$ 
Summing over all choices of markings gives
$$
\sum_{J\subseteq\{1, \ldots, b\}}(-2)^{|J|}
=\sum_{k=0}^{b}\binom bk(-2)^k
=(1-2)^b=(-1)^b,
$$
Each marking corresponds uniquely, in the graph $G$, to a term for the original matching $M'$ in the expansion in \eqref{eq:2.1}. 
Therefore,
\begin{equation}\label{eq:4.11}
P_G(\mathbf{q})=\sum_{G'}(-1)^{b}
\prod_{v_i\in R(G')}t_{i}',
\end{equation}
Here we no longer distinguish whether the pairs of opposite arrows in $G'$ arise from matching edges,
 and $R(G')$ denotes the set of all roots of $G'$.

Now retain only the edges of $G'$ carrying arrows (deleting all edges without arrows), 
and then remove all arrows.
 In the spanning subgraph of $G'$ obtained by this procedure, the number of edges in each connected component $C$ is at most the number of outgoing arrows,
 so $|E(C)|\le|V(C)|$.
 Thus this spanning subgraph is a pseudo-forest.
 Fix a pseudo-forest. We now sum all contributions on its connected components.

If a connected component $C$ is a tree $T$, let $s'$ be its number of vertices, $r'$ its number of roots, and $b_T$ its number of pairs of opposite arrows.
 Counting the total number of outgoing arrows in two ways gives $s'-r'$ and $|E(T)|+b_T=s'-1+b_T$, respectively. Hence $r'=1-b_T$.
 Thus $r'$ can only be $1$ or $0$:
  when $r'=1$, we have $b_T=0$. After choosing the unique root $v_i$ arbitrarily, all other edges must point toward this root along the unique paths in the tree,
  giving a total weight of $\sum_{v_i\in V(T)}t_{i}'$;
  when $r'=0$, we have $b_T=1$. After choosing the unique edge carrying opposite arrows arbitrarily, all other arrows must point toward this edge. Each choice has weight $-1$, giving a total weight of $-|E(T)|$. Therefore, the total contribution of a tree component is
$$
\sum_{v_i\in V(T)}t_{i}'-|E(T)|
=1+\sum_{v_i\in V(T)}a_{i}.
$$

If a connected component $C$ is unicyclic, let $r''$ be its number of roots and $b_C$ its number of pairs of opposite arrows.
 The same arrow count gives $|V(C)|-r''=|E(C)|+b_C=|V(C)|+b_C$, so $r''=b_C=0$.
  Its unique cycle must be oriented consistently, and all arrows outside the cycle must point toward it.
  Since the cycle has two possible orientations, the total contribution is $2$.

Since the choices for different components are independent,
 combining \eqref{eq:4.11} with a sum over all spanning pseudo-forests and taking the product of the contributions of their components
  gives \eqref{eq:2.2}, completing the proof.
\end{proof}
\begin{remark}
When all $a_i$ vanish, Lemma \ref{lem:2.1} reduces to the pseudo-forest expansion of the polynomial $R_G(z)$ at $z=2$ given by Bencs and Csikv\'ari \cite{Bencs2022,Bencs2023}.

 To illustrate Lemma~\ref{lem:2.1}, consider the graph $G = G(1,0,1)$ (defined in the introduction), which is a triangle with one pendant edge, and set $a_i = 0$ for every vertex. All $16$ spanning subgraphs of $G$ are pseudo-forests. Each of the two spanning subgraphs containing all three edges of the triangle has exactly one unicyclic component, while the remaining $14$ are forests.
  Thus, Eq.~\eqref{eq:2.2} yields $P_G(\mathbf{q}) = 14 + 2 \cdot 2 = 18$. 
In Eq.~\eqref{eq:2.1}, the empty matching contributes $72$, the four single-edge matchings together contribute $-2(6 + 6 + 9 + 8)$, and the unique two-edge matching contributes $4$, yielding the same value.
\end{remark}

\begin{lemma}\label{lem:2.3}
Let $G$ be a graph on $n$ vertices, and let $H$ be the graph obtained from $G$ by deleting a vertex $v (v\in V(G))$. 
 For each $v_i\in V(H)$,
 let $q_{i}=d_{i}+1$. For a vertex set $S\subseteq V(H)$, define
$$
  R_{i}(S)=\frac{P(S)}{P(S\setminus\{v_i\})} \quad (v_i\in S), \qquad
  b_{i}(S)=q_{i}-2\sum_{v_j\in N_{H[S]}(v_i)}\frac1{q_{j}},
$$
where 
$
  P(S)=P_{H[S]}(\mathbf{q}|_S)$, $
  \mathbf{q}|_S=(q_i)_{v_i\in S}$, 
and $\mathbf{q}|_S$ inherits the weights from $G$.
 Then
\begin{equation}\label{eq:2.3}
F(G)=\sum_{\varphi\in\mathcal{F}(H)}\prod_{T\in\mathcal{T}(\varphi)}\left(1+|N_G(v)\cap V(T)|\right)\le P_H(\mathbf{q}),
\end{equation}
\begin{equation}\label{eq:2.4}
0<R_{i}(S)=q_{i}-2\sum_{v_j\in N_{H[S]}(v_i)}\frac1{R_{j}(S\setminus\{v_i\})}
\le q_{i},
\end{equation}
\begin{equation}\label{eq:2.5}
R_{i}(S)\le b_{i}(S),
\qquad b_{i}(S)\ge 1,
\end{equation}
Equality in \eqref{eq:2.3} holds if and only if $H$ is a forest.
\end{lemma}
\begin{proof}
For $H$, fix a spanning forest $\varphi$. We now add the vertex $v$ back to obtain spanning forests of $G$.
 In $\varphi$, fix a tree component $T_i$, and let the number of its vertices adjacent to $v$ be
$$
a_i'=|N_G(v)\cap V(T_i)|.
$$
There are two types of extensions for this tree component:
the first is not to connect it to $v$, for which there is only $1$ possibility;
the second is to connect one vertex of $T_i$ to $v$, for which there are $a_i'$ possibilities.
If $T_i$ has $2$ or more vertices connected to $v$, a cycle is formed, contrary to the definition of a spanning forest. 
The two choices are mutually exclusive, so there are 
$
1+a_i'=1+|N_G(v)\cap V(T_i)|
$ choices in total. 
Thus, for the fixed $\varphi$, the number of spanning forests of $G$ obtained is
$$
\prod_{T\in\mathcal{T}(\varphi)}
\left(1+|N_G(v)\cap V(T)|\right).
$$
For any spanning forest of $G$, deleting the vertex $v$ and its incident edges uniquely determines a spanning forest of $H$, denoted by $\varphi$. Conversely, given $\varphi$ and the choices above for each of its tree components,
 a spanning forest of $G$ is uniquely determined. Hence
 there is a bijection between the spanning forests of $G$ and all $\varphi$ together with the above choices.
 Summing over $\mathcal{F}(H)$ gives
 $$
 F(G)=\sum_{\varphi\in\mathcal{F}(H)}\prod_{T\in\mathcal{T}(\varphi)}\left(1+|N_G(v)\cap V(T)|\right).
 $$
 
The above analysis shows that, for each $v_i\in V(H)$,
$$
q_{i}=d_H(v_i)+1+\varepsilon_{i}, \qquad \varepsilon_{i}=1_{\{v{v_i}\in E(G)\}}\ge0.
$$
Applying Lemma \ref{lem:2.1} gives 
$a_{i}=d_{i}-d_H(v_i)$, and hence
 $1+\sum_{v_i\in V(T)}a_{i}=1+|N_G(v)\cap V(T)|$.
  Thus the sum in \eqref{eq:2.3} consists of exactly all the forest terms in \eqref{eq:2.2}. 
  Therefore, $F(G)\le P_H(\mathbf{q})$, with equality in \eqref{eq:2.3} if and only if $H$ is a forest.

For any vertex set $S\subseteq V(H)$ and $v_i\in S$,
 let $a_{i}=q_{i}-d_{H[S]}({v_i})-1=d_{i}-d_{H[S]}(v_i)\ge0$.
 By Lemma \ref{lem:2.1}, both $P(S)$ and $P(S\setminus\{v_i\})$ are greater than $0$,
 so $R_{i}(S)>0$. 
 We divide the matchings of $H[S]$ into two classes according to whether they cover the vertex $v_i$,
  and use \eqref{eq:2.1} to derive a recurrence for $P(S)$.
 First, consider matchings that do not cover the vertex $v_i$. 
In this case, $v_i$ contributes a weight of $q_{i}$ to the weighted sum over matchings.
The matching on the remaining vertices can be any matching of $H[S\setminus\{v_i\}]$.
Thus the total contribution of this class of matchings is $q_{i}P(S\setminus\{v_i\})$.
Next, consider matchings that cover the vertex $v_i$.
 In this case, exactly one matching edge is incident with $v_i$.
Let this edge be $v_iv_j$; it contributes a weight of $-2$.
Since both $v_i$ and $v_j$ are covered by this matching edge, the remaining matching edges can only be chosen from the graph
 $H[S\setminus\{v_i, v_j\}]$.
 Therefore, for a fixed $v_j$, the total contribution of these matchings is
 $-2P(S\setminus\{v_i, v_j\})$. 
For $v_i$, summing over all its neighbors $v_j$ gives
$$
-2\sum_{v_j\in N_{H[S]}(v_i)}
P(S\setminus\{v_i, v_j\}).
$$
Combining the two cases above, we obtain
\begin{equation}\label{eq:3.10}
P(S)=q_{i}P(S\setminus\{v_i\})-2\sum_{v_j\in N_{H[S]}(v_i)}P(S\setminus\{v_i, v_j\}).
\end{equation}
Dividing both sides of \eqref{eq:3.10} by $P(S\setminus\{v_i\})$ gives \eqref{eq:2.4}.
 In $S\setminus\{v_i\}$, take the vertex $v_j$, apply \eqref{eq:2.4}, and take reciprocals to obtain
$$
\frac1{R_{j}(S\setminus\{v_i\})}\ge\frac1{q_{j}},
$$
Substituting this into \eqref{eq:2.4} gives $R_{i}(S)\le b_{i}(S)$.
 For the vertex $v_i$ and any of its neighbors $v_j$ in the sum,
 since $d_j\ge1$, we have $q_{j}\ge2$. Hence
 
$$
b_{i}(S)
=q_{i}
-\sum_{v_j\in N_{H[S]}(v_i)}
\frac{2}{q_{j}}
\ge q_{i}
-\sum_{v_j\in N_{H[S]}(v_i)}1
=d_{i}+1-d_{H[S]}(v_i)
\ge1,
$$
This completes the proof.
\end{proof}
We now prove Theorem \ref{the:2}.
\begin{proofthm}{Proof of Theorem \ref{the:2}}
When $n=1$, $G$ is edgeless, $F(G)=1$, and the right-hand side of \eqref{eq:1.3} is also $1$.
 Thus Theorem \ref{the:2} clearly holds. Henceforth, assume $n\ge2$.

 For $2\le i\le n$, let $q_{i}=d_i+1$,
 $S_i=\{v_i, \ldots, v_n\}$, $S_{n+1}=\varnothing$.
  Let $H'=G-v_1$.
  In Lemma \ref{lem:2.3}, take $v=v_1$ and $H=H'$.
   We retain its notation $P(S)$ and $R_i(S)$, 
  with $q_i$ always equal to $d_i+1$.
   By \eqref{eq:2.3}, we obtain
\begin{equation}\label{eq:3.1}
F(G)\le P_{H'}(\mathbf{q})=\prod_{i=2}^n\frac{P(S_i)}{P(S_{i+1})}=\prod_{i=2}^n R_{i}(S_i).
\end{equation}
Write
$$
N_i=\{j:j>i \ \text{and} \ v_iv_j\in E(G)\}, \qquad b_i=q_i-2\sum_{j\in N_i}\frac1{q_j}.
$$
By \eqref{eq:2.5}, we have
\begin{equation}\label{eq:3.2}
0<R_{i}(S_i)\le b_i, \qquad b_i\ge1.
\end{equation}
 The neighbors of $v_i$ among $v_1, \ldots, v_{i-1}$ number at most $i-1$. Hence
\begin{equation*}
|N_i|\ge\max\{d_i-(i-1), 0\}=k_i.
\end{equation*}
Since $q_2\ge q_3\ge\cdots\ge q_n>0$, we have
$$
\frac1{q_2}\le\frac1{q_3}\le\cdots\le\frac1{q_n}.
$$
 If $k_i>0$, arrange the elements of $N_i$ in increasing order 
 as $j_1<\cdots<j_{|N_i|}$. Then $j_\ell\ge i+\ell$, so
\begin{equation}\label{eq:27}
\sum_{j\in N_i}\frac1{q_j}\ge\sum_{\ell=1}^{k_i}\frac1{q_{j_\ell}}\ge\sum_{\ell=1}^{k_i}\frac1{q_{i+\ell}}=h_i.
\end{equation}
 When $k_i=0$, the corresponding $h_i=0$. Hence, for every $k_i\geq0$, we have
\begin{equation}\label{eq:1}
\sum_{j\in N_i}\frac1{q_j}\ge h_i.
\end{equation}
 Therefore, for every $i\ge2$, we have
\begin{equation}\label{eq:3.4}
1\le b_i\le q_i-2h_i=d_i+1-2h_i.
\end{equation}
All factors in \eqref{eq:3.1}, \eqref{eq:3.2}, and \eqref{eq:3.4} are positive, so the inequalities can be multiplied term by term to give
$$
F(G)\le\prod_{i=2}^n b_i\le\prod_{i=2}^n(d_i+1-2h_i),
$$
This proves \eqref{eq:1.3}.
We now discuss the case of equality in Theorem \ref{the:2}.

First, we prove necessity. Suppose that equality holds in \eqref{eq:1.3} and that $G$ is not edgeless. Then $d_2\ge1$.
By \eqref{eq:3.1}, \eqref{eq:3.2}, and \eqref{eq:3.4}, for every $2\le i\le n$, we have
$$
R_{i}(S_i)=b_i, \qquad \sum_{j\in N_i}\frac1{q_j}=h_i.
$$
Equality in \eqref{eq:1} requires $|N_i|=k_i$. 
 Thus, when $k_i=0$, we have $N_i=\varnothing$.
When $k_i>0$, term-by-term comparison gives
\begin{equation}\label{eq:30}
d_{j_\ell}=d_{i+\ell}, \qquad 1\le\ell\le k_i.
\end{equation}
On the other hand, \eqref{eq:2.4} gives
$$
b_i-R_{i}(S_i)
=2\sum_{j\in N_i}\left(\frac1{R_{j}(S_{i+1})}-\frac1{q_j}\right)=0,
$$
Thus, for every $j\in N_i$, we have
$R_{j}(S_{i+1})=q_j$. For $R_j(S_{i+1})$, applying \eqref{eq:2.4} again gives
$$
N_{G[S_{i+1}]}(v_j)=\varnothing, \qquad j\in N_i.
$$

Since $|N_2|=k_2=d_2-1$, we have $v_1v_2\in E(G)$.
 We consider two cases for $d_2$:
if $d_2=1$, then for all $i\ge3$, we have $d_i\le1$;
if $d_2\ge2$, then $k_2>0$, and for each $j\in N_2$, the corresponding vertex $v_j$ is isolated in $G[S_3]$,
so its neighbors can only belong to $\{v_1,v_2\}$, and hence $d_j\le2$.
In \eqref{eq:30}, taking $i=2$ and $\ell=1$,
 and using the equality conditions in \eqref{eq:2.4}, \eqref{eq:3.2}, and \eqref{eq:27}, 
 we obtain $d_3=d_{j_1}\le2$
 (where $j_1$ is the smallest index in $N_2$).
Thus, for all $i\ge3$, we have
$$
|N_i|=  k_i=(d_i-i+1)_+=0,
$$
Hence $G-\{v_1,v_2\}$ is edgeless.
Among the vertices of $V(G)\setminus\{v_1,v_2\}$, let $t_1,t_2,t_3$ denote the numbers of vertices adjacent only to $v_1$, only to $v_2$, and to both $v_1,v_2$, respectively.
All remaining vertices are isolated.
If $t_3\ge1$, the graph obtained after deleting the isolated vertices is precisely $G(t_1,t_2,t_3)$.
If $t_3=0$ and $t_1,t_2\ge1$, the resulting graph is the double star $D_{t_1,t_2}$.
If $t_3=0$ and at least one of $t_1,t_2$ is $0$, the resulting graph is the star $K_{1,t_1+t_2+1}$.
This proves necessity.

Next, we prove sufficiency. If $G$ is edgeless, then $F(G)=1$, and every factor on the right-hand side of \eqref{eq:1.3} is $1$, so equality holds.

 For the star $K_{1,t}$, 
 the number of spanning forests is $F(K_{1,t})=2^{t}=2^{n-1}$.
 Since $d_i=1(2\le i\le n)$
  and $k_i=h_i=0(2\le i\le n)$, 
 we have $\prod_{i=2}^n(d_i+1-2h_i)=2^{n-1}$,
 that is, $$F(K_{1,t})=\prod_{i=2}^n(d_i+1-2h_i).$$
 
 For the double star $D_{p,q}(p\geq q)$,
  the number of spanning forests is $F(D_{p,q})=2^{p+q+1}=2^{n-1}$.
  Since $$d_2=q+1, \qquad d_i=1(3\le i\le n), \qquad k_2=q,  \qquad h_2=\frac{q}{2}, 
   \qquad k_i=h_i=0(3\le i\le n),$$
  we have $\prod_{i=2}^n(d_i+1-2h_i)=2^{n-1}$,
 that is, $$F(D_{p,q})=\prod_{i=2}^n(d_i+1-2h_i).$$
 
Without loss of generality, assume that $t_1\ge t_2$. Upon deleting the $t_1+t_2$ pendant vertices, the graph $G(t_1,t_2,t_3)$ reduces to $G(0,0,t_3)$, namely the book graph consisting of $t_3$ triangles sharing a common edge. By \cite{Shteiner2025},
$$
F(G(0,0,t_3))=2\cdot 3^{t_3}+t_3\,3^{t_3-1}.
$$
Since every pendant edge is a bridge and hence lies on no cycle, the choice of including or omitting each pendant edge is independent of the choices made for the remaining pendant edges. Consequently,
$$
F(G(t_1,t_2,t_3))
=2^{t_1+t_2}\bigl(2\cdot 3^{t_3}+t_3\,3^{t_3-1}\bigr)
=2^{t_1+t_2}3^{t_3}\left(2+\frac{t_3}{3}\right).
$$
 Since
$$
 d_2=t_2+t_3+1,\qquad d_i=2(3\leq i \leq t_3+2), \qquad d_i=1(t_3+3\leq i \leq n),
$$
we have
$$
k_2=t_2+t_3, \qquad h_2=\frac{t_3}{3}+\frac{t_2}{2}, \qquad
k_i=h_i=0\quad(3\le i\le n),
$$
that is,
$$
\prod_{i=2}^n(d_i+1-2h_i)
=\left(t_2+t_3+2-\frac{2t_3}{3}-t_2\right)3^{t_3}2^{t_1+t_2}
=2^{t_1+t_2}3^{t_3}\left(2+\frac{t_3}{3}\right)=F(G).
$$
Adding isolated vertices does not change the value of either side of \eqref{eq:1.3}. Hence all the graph classes above attain equality,
 proving sufficiency.

\end{proofthm}

\section{Further discussion}
Sections 2 and 4 give two different upper bounds and their corresponding extremal graphs, respectively.
 A natural question arises: for which families of connected graphs with  $n$ vertices and 
$m$  edges is one of the bounds smaller than the other?
 
For convenience of comparison, for any connected graph with $n\ge2$ vertices and $m$ edges, denoted by $G$, write
$$
U_d(G)=\prod_{i=2}^n(d_i+1-2h_i), \qquad
U_{\mu}(G)=\tau(G)\bigl((1+\mu^{-1})^n-\mu^{-n}\bigr).
$$
We use several graph classes below to show that neither of these two upper bounds is uniformly better than the other.

For $t\ge2$, the star $K_{1,t}$ has algebraic connectivity $1$, and hence
$$
U_d(K_{1,t})=F(K_{1,t})=2^{n-1}<2^n-1=
U_{\mu}(K_{1,t}).
$$

For the double star $D_{p,q}$ and, when $n\ge4$, the graph $G(t_1,t_2,t_3)$,
equality holds in Theorem \ref{the:2}, but not in Theorem \ref{cor:1.1}, and therefore
$$
U_d(D_{p,q})=F(D_{p,q})<U_{\mu}(D_{p,q}),
$$
$$
U_d(G(t_1,t_2,t_3))=F(G(t_1,t_2,t_3))<U_{\mu}(G(t_1,t_2,t_3)).
$$

For $s\ge2$, a calculation for the complete bipartite graph
$K_{s,s}$ gives
$$
U_d(K_{s,s})
=(s+1)^s\prod_{i=2}^{s}\left(s+1-\frac{2(s-i+1)}{s+1}\right)
=(s+1)^{2s-1}\prod_{j=1}^{s-1}
\left(1-\frac{2j}{(s+1)^2}\right).
$$
For $0\le x_j\le1$, we have
 $\prod_j(1-x_j)\ge1-\sum_jx_j$, 
and applying this to the expression above gives
$$
\prod_{j=1}^{s-1}\left(1-\frac{2j}{(s+1)^2}\right)
\ge1-\frac{s(s-1)}{(s+1)^2}
=\frac{3s+1}{(s+1)^2},
$$
Thus,
 $$
(3s+1)(s+1)^{2s-3}
\le U_d(K_{s,s}).$$
The Laplacian eigenvalues of the complete bipartite graph $K_{s,s}$ are
$$
2s,\ \underbrace{s,\ldots,s}_{2s-2\text{ copies}},\ 0,
$$
so its algebraic connectivity is $s$.
From $n\tau(G)=\prod_{i=1}^{n-1}\lambda_i$, we obtain
 $\tau(K_{s,s})=s^{2s-2}$, 
and therefore
$$
U_{\mu}(K_{s,s})=s^{2s-2}
\left(\left(1+\frac1s\right)^{2s}-\frac1{s^{2s}}\right)
=\frac{(s+1)^{2s}-1}{s^2}.
$$
Moreover, since
$$
s^2(3s+1)-(s+1)^3
=(s-2)(2s^2+2s+1)+1>0,
$$
we have
$$
U_{\mu}(K_{s,s})
<\frac{(s+1)^{2s}}{s^2}
<(3s+1)(s+1)^{2s-3}
\le U_d(K_{s,s}).
$$
Thus, for this family of complete bipartite graphs, the upper bound given by Theorem \ref{cor:1.1} is strictly better than that given by Theorem \ref{the:2}.

The above discussion shows that neither upper bound is uniformly tighter than the other. Consequently, a complete characterization of graph classes in terms of the relative sizes of these two upper bounds remains an open problem.
 
\noindent{\bf Conflict of Interest Statement}

The authors declare that they have no conflicts of interest.


\begin{thebibliography}{99}

\bibitem{Aldous2007}
D. Aldous, R. Lyons, Processes on unimodular random networks,  {\em Electron. J. Probab. } 12 (2007) 1454--1508.

\bibitem{Annan1994}
J.D. Annan, A randomised approximation algorithm for counting the number of forests in dense graphs,  {\em Combin. Probab. Comput. } 3 (1994) 273--283.

\bibitem{Bencs2022}
F. Bencs, P. Csikv\'ari, Evaluations of Tutte polynomials of regular graphs, {\em J. Combin. Theory Ser. B} 157 (2022) 500--523.
\bibitem{Bencs2023}
F. Bencs, P. Csikv\'ari, Upper bound for the number of spanning forests of regular graphs,  {\em European J. Combin. } 110 (2023) 103677.

\bibitem{Borbenyi2021}
M. Borb\'enyi, P. Csikv\'ari, H. Luo, On the number of forests and connected spanning subgraphs,  {\em Graphs Combin. } 37 (2021) 2655--2678.

\bibitem{Caracciolo2004}
S. Caracciolo, J.L. Jacobsen, H. Saleur, A.D. Sokal, A. Sportiello, Fermionic field theory for trees and forests,  {\em Phys. Rev. Lett. } 93 (2004) 080601.

\bibitem{Chebotarev1997}
P. Chebotarev, E.V. Shamis, The matrix-forest theorem and measuring relations in small social groups,  {\em Autom. Remote Control} 58 (1997) 1505--1514.
\bibitem{Diestel2017}
R. Diestel, Graph theory, fifth ed. , in: Graduate Texts in Mathematics, vol. 173, Springer, Berlin, 2017, p. xviii+428.
\bibitem{Grimmett2004}
G.R. Grimmett, S.N. Winkler, Negative association in uniform forests and connected graphs,  {\em Random Struct. Algorithms} 24 (2004) 444--460.
\bibitem{Horn2013}
R.A. Horn, C.R. Johnson,
Matrix Analysis, second ed., Cambridge University Press, Cambridge, 2013.
\bibitem{Jaeger1990}
F. Jaeger, D.L. Vertigan, D.J.A. Welsh, On the computational complexity of the Jones and Tutte polynomials,  {\em Math. Proc. Cambridge Philos. Soc. } 108 (1990) 35--53.

\bibitem{Kelmans1974}
A.K. Kelmans, V.M. Chelnokov, A certain polynomial of a graph and graphs with an extremal number of trees, {\em J. Combin. Theory Ser. B} 16 (1974) 197--214.
\bibitem{Liu1981}
C.J. Liu, Y. Chow, Enumeration of forests in a graph,  {\em Proc. Amer. Math. Soc. } 83 (1981) 659--662.
\bibitem{Luo2026}
J. Luo, T. Wu, Y.-N. Yeh, Counting spanning trees of multipartite graphs and those containing a fixed spanning forest,  {\em Discrete Appl. Math.} 395 (2026) 317--336.
\bibitem{Shteiner2025}
S. Shteiner, P. Shteyner, Comparing the numbers of subforests and
subgraph-degree-tuples, arXiv:2510.26936v1 (2025).
\bibitem{Stark2013}
D. Stark, The asymptotic number of spanning forests of complete
bipartite labelled graphs, {\em Discrete Math. } 313 (2013) 1256--1261.
\bibitem{Teranishi2005}
Y. Teranishi, The number of spanning forests of a graph, {\em Discrete Math. } 290 (2005) 259--267.

\bibitem{Thomassen2010}
C. Thomassen, Spanning trees and orientations of graphs, {\em J. Comb.} 1 (2010) 101--111.

\bibitem{Wu2026}
T. Wu, S. Lu, X. Dong, A sharp upper bound on the number of spanning forests of regular graphs, preprint, arXiv:2610.05722v1 (2026).
\bibitem{Xu2026}
S. Xu, K. Xu, On asymptotic values for the minimum number of spanning forests in simple regular graphs, preprint, arXiv:2605.19707v2 (2026).
\end{thebibliography}
\end{document}